\documentclass[a4paper,11pt]{amsart}
\usepackage[cp1250]{inputenc}
\usepackage{amssymb}
\usepackage{amsfonts}
\usepackage{exscale}
\usepackage{amsmath}
\usepackage[T1]{fontenc}

\newtheorem{defin}{Definition}
\newtheorem{twr}[defin]{Theorem}
\newtheorem{lem}[defin]{Lemma}

\newtheorem{prop}[defin]{Proposition}

\def \N {\mathbb{N}}
\def \ndon {\N^{n}}

\def \field {\mathbb{K}}
\def \dowod {\noindent {\sc Proof. }}
\def \remark {\noindent {\sc Remark. }}

\def \M {\mathcal{M}}
\def \sgn {\operatorname{sgn}}
\def \impl {\Longrightarrow}
\def \degred {\operatorname{degred}}
\def \S {\mathcal{S}}
\def \kxn {\field[x_{1},\dots,x_{n}]}
\def \grF {\mathbf{F}}
\def \remarkspace {\vspace{0.3cm}}

\begin{document}

\title[Expected term bases...]{Expected term bases for generic \\ multivariate Hermite interpolation}

\author{Marcin Dumnicki}

\dedicatory{
Institut of Mathematics, Jagiellonian University, \\
Reymonta 4, 30-059 Krak\'ow, Poland \\
E-mail address: Marcin.Dumnicki@im.uj.edu.pl
}


\begin{abstract}
The main goal of the paper is to find an effective estimation for
the minimal number of generic points in $\field^{2}$ for which the basis
for Hermite interpolation consists of the first $\ell$ terms (with
respect to total degree ordering). As a result we prove that the
space of plane curves of degree $d$ having generic singularities
of multiplicity $\leq m$ has the expected dimension if the number
of low order singularities (of multiplicity $k \leq 12$) is greater
then some $r(m,k)$. Additionally, the upper bounds for $r(m,k)$ are given.
\end{abstract}

\maketitle

\section{Introduction.}

We denote by $\N$ the set of nonnegative integers, by $\field$ a field
of characteristic zero. Let $I_{n} := \{ 1,2, \dots, n \}$.
We will use the natural one-to-one correspondence
between monomials $x^\alpha \in \field[x_{1}, \dots, x_{n}]$
and multiindices $\alpha \in \N^{n}$. For any two multiindices
$\alpha, \beta$ we will write $\beta \leq \alpha$ if
$\alpha - \beta$ has only nonnegative entries.

By a {\it Ferrers diagram $F$} we understand a finite
subset $F \subset \N^{n}$
such that if $\alpha \in F$, $\beta \leq \alpha$ then $\beta \in F$.

Let $\mathbf{F} = \{F_{j}\}_{j=1}^{r}$ be a finite sequence of Ferrers diagrams,
let $\mathbf{P} = \{p_{j}\}_{j=1}^{r}$ be the sequence of parwise
different points in $\field^{n}$.
The interpolation ideal assigned to $\mathbf{F}$ and $\mathbf{P}$ is the ideal
$$I = \big\{ f \in \kxn : \frac{\partial^{|\alpha|}f}{\partial x^{\alpha}}(p_{i}) = 0,
\ \alpha \in F_{i}, i=1,\dots,r \big\}.$$
Let us introduce the multivariate Hermite interpolation problem, that
is the problem of finding a basis $B = B(\mathbf{F},\mathbf{P})$ of 
$\kxn / I$ as a vector space over $\field$. The classical approach
is to compute the Gr\"obner basis of $I$ with respect to an admissible
ordering (cf. \cite{WF}). This method gives a minimal
basis (with respect to the chosen admissible ordering) of the quotient
space. However, due to time complexity, it is not very practical.

Consider the sequence of Ferrers diagrams $\mathbf{F}$ and an admissible
ordering. The basis $B$ depends on
the sequence of points $\mathbf{P} \in (\field^{n})^{r}$, 
but there exists one special
basis (called the {\it generic basis}) which is the same for
almost all $\mathbf{P}$, that is for $\mathbf{P}$ in a Zariski open,
dense subset of $(\field^{n})^{r}$. The problem of finding the 
interpolation basis
(generic or not) for the lexicographical ordering without
using the Buchberger's algorithm was solved in \cite{CM} (the non-generic case)
and \cite{ASTW} (the generic case).

For a total degree ordering the methods of finding the interpolation
basis without the Buchberger algorithm are not known.
An important question that arises here is:

How can we characterize the sequences of Ferrers diagrams
for which the generic basis $B$ of interpolation
is contained in the set $\{\alpha : |\alpha| \leq d \}$?
If we assume that all Ferrers
diagrams are of the form $\{\alpha : |\alpha| \leq m\}$ then
this problem is closely
related to the problem of finding the actual dimension
of the space of hypersurfaces (in $\field^{n}$)
of degree $d$ having generic
singularities of multiplicity $m$ (homogeneous generic singularities problem)
or up to multiplicity $m$ (inhomogeneous generic singularities problem).

The last problem was solved
by J. Alexander and A. Hirschowitz
(\cite{AH1}, \cite{AH2}) who showed that for the number of
singularities large enough this dimension is the expected dimension,
however they do not give a bound for the number of singularities needed.
For some cases the problem was studied in many other papers.
The homogeneous case for $n = 2$, $m \leq 12$
is completed in \cite{CMir}, the inhomogeneous case for $n = 2$, $m \leq 4$
in \cite{Mig}.
A more computational approach to this problem can be found
in \cite{S} and \cite{MS}.

We present an effective criterion for the sequence $\mathbf{F}$
to have the desired form of the basis $B$. As a result
we present new proofs
for the inhomogeneous generic singularities problem for $m \leq 12$ together
with the bound for the number of singularities needed.
Moreover, for arbitrary $m$ we give the bound for sufficient number of singularities
of multiplicity $k \leq 12$:

\begin{twr}
\label{curvethm}
Let $\Gamma_{d,p_{0},\dots,p_{m}}$ denote the space of
all plane curves of degree $d$ passing through $p_{0}$ generic points
and having $p_{j}$ generic singularities of order $j$, for
$j = 1, \dots, m$. Let $0 \leq k \leq 12$, $k \leq m$.
There exists $r(m, k)$, 
$$ r(m,k) \leq \max \left\{ 6(m+1), \frac{4(m+1)(2m+1)}{(k+1)(k+2)} \right\}$$
such that if $p_{k} > r(m, k)$ then 
$\Gamma_{d,p_{0},\dots,p_{m}}$ has the expected dimension
(as a vector space over base field) equal to
$$\dim \Gamma_{d,p_{0},\dots,p_{m}} = \max \left\{ 0, \frac{(d+1)(d+2)}{2} -
\sum_{k=0}^{m} p_{k}\frac{(k+1)(k+2)}{2} \right\}.$$
\end{twr}

We discuss the method of finding such bounds, and present the strict
values for $k,m \leq 7$. Our method is a new one, we do not refer to the methods
used in other papers.

In sections 2--4 we introduce the methods
and prove lemmas used in section 5, which is the main section for this paper.
An example of using our method
for finding the basis $B$ for arbitrary sequence of Ferrers diagrams
appears in section 6. 

\section{Generically correct problems.}

For any monomial $x^{\alpha} \in \field[x_{1},\dots,x_{n}]$, $
\alpha \in \ndon$,
a multiindex $\beta \in \ndon$ and
a point $a = (a_{1}, \dots, a_{n}) \in \field^{n}$
we define
$$\varphi(x^{\alpha}, \beta, a) = \left\{
\begin{array}{ll} \frac{\alpha_{1}! \cdot \cdots \cdot \alpha_{n}!}
{(\alpha_{1}-\beta_{1})! \cdot \cdots \cdot (\alpha_{n}-\beta_{n})!}
a_{1}^{\alpha_{1}-\beta_{1}} \cdot
\cdots \cdot a_{n}^{\alpha_{n}-\beta_{n}}, & \textrm{ if } \beta \leq \alpha, \\
0, & \textrm{ otherwise. } \end{array} \right.$$
$\varphi(x^{\alpha}, \beta, a)$ is just a partial derivative of $x^{\alpha}$ with
respect to $\beta$ taken at the point $a$.

Let $p_{1}, \dots, p_{r}$ be $r$ distinct
points (nodes) in $\field^{n}$. Let $\grF = \{F_{i}\}_{i=1}^{r}$.
Define a set of conditions
$$C_{\grF} = \{ (p, \beta) \in \field^{n} \times \ndon
 : \exists_{i \in I_{r}} \ p = p_{i}, \ \beta \in F_{i} \}.$$
The cardinality
of $C_{\grF}$ (denoted by $c$) is equal to the sum $\sum_{i=1}^{r} \# F_{i}$.
Let us assume that the set of monomials $B \subset \field[X_{1},\dots,X_{n}]$
of cardinality $c$ is given.
We can order sets $C_{\grF} = \{ (p,\beta)_{1}, \dots, (p, \beta)_{c} \}$
and $B = \{ x^{\alpha_{1}}, \dots, x^{\alpha_{c}} \}$ and define the matrix
$$M = \big[ a_{k,\ell} \big]_{k=1,\dots,c}^{\ell=1,\dots,c},$$
where $a_{k,\ell} = \varphi(x^{\alpha_{\ell}}, \beta, p)$, $
(p, \beta) = (p, \beta)_{k}$.
We say that {\it the interpolation problem for the sequence of Ferrers
diagrams $\grF$
and the set of monomials $B$ is correct (shortly $(\grF, B)$ is correct)}
if $\det M \neq 0$.
Of course the correctness of the interpolation problem does not depent
on ordering of $\grF$ and $B$.

Let us make the following observation. The interpolation problem is correct
if and only if the following is true:

For any set of values (of cardinality $c$)
we can find a polynomial $P$ in the linear space
spanned by $B$ over $\field$ having prescribed values and
derivatives in each node. The matrix $M$ is just the matrix of the linear
equation solving this problem, and $B$ is the basis of the quotient
space $\kxn / I$, where $I$ is the interpolation ideal.

The determinant of the matrix $M$ can be considered as a polynomial
of $nr$ coordinates of nodes, say $\det M \in \field[p_{1}^{1}, \dots,
p_{1}^{n}, p_{2}^{1}, \dots, p_{2}^{n}, \dots, p_{r}^{1}, \dots, p_{r}^{n}]$.
We say that {\it the interpolation problem $(\grF, B)$ is generically correct}
if $\det M$ is a nonzero polynomial.

Observe that $(\grF, B)$ is generically correct if and only if
there exists a set of nodes $\mathbf{P}$ for which $(\grF, B)$
is correct, and if and only if
it is correct for the set of points from a Zariski open, 
dense subset of $(\field^{n})^{r}$. Hence $B$ is a generic basis for
interpolation.

Let $F$ be a Ferrers diagram, $B$ a set of monomials,
$B' = \{ x^{\alpha_{1}}, \dots, x^{\alpha_{k}} \} \subset B$,
$\# B' = \# F$. We say that {\it $B'$ is exceptional in $B$ with respect to $F$}
if the following conditions are fulfilled:\\
1. For any $P = \{ x^{\beta_{1}}, \dots, x^{\beta_{k}} \} \subset B$, $P \neq B'$
such that $\prod_{i=1}^{k} x^{\alpha_{i}} = \prod_{i=1}^{k} y^{\beta_{i}}$
the problem $(\{F\},P)$ is
not generically correct,\\
2. The problem $(\{F\},B')$ is generically correct.

\begin{twr}
\label{maintred}
Let $(\grF, B)$ be a generically correct interpolation problem,
let $P$ be a set of monomials, let $F$ be a Ferrers diagram.
Denote $\grF' = \grF \cup \{F\}$ (this is not the sum of sets, but adding
an element to the sequence), $B' = B \cup P$. If
$B \cap P = \varnothing$, $\#F = \#P$ and $P$ is exceptional in $B'$
with respect to $F$, then the problem $(\grF', B')$ is generically
correct.
\end{twr}

In the proof we will use the following Lemma:

\begin{lem}[generalized Laplace rule]
\label{leibnitzlem}
Let $M \in \M(n,n;\field)$ be a square matrix, let $k \in \N$, $1 \leq k < n$.
Denote
$$M = \left[ \begin{array}{c} M_{1} \\ M_{2} \end{array} \right],$$
where $M_{1} \in \M(k,n;\field)$ and $M_{2} \in \M(n-k,n;\field)$.
Let
$$\S = \{ (a_{1}, \dots, a_{k}) : a_{i} \in I_{n}, a_{1} < a_{2} < \dots < a_{k} \}$$
be the set of all possible chosing of $k$ columns (without order).
Let $m_{1}^{S}$ be the minor of $M_{1}$ determined by $S \in \S$,
$m_{2}^{S}$ be the $n-k \times n-k$ minor of $M_{2}$ determined by
$I_{n} \setminus S$. Then
$$\det M = \sum_{S \in \S} \sgn(s(S))m_{1}^{S}m_{2}^{S},$$
where $\sgn(s(S)) = \pm 1$, $s(S)$ being a permutation
$$s(S) = \left( \begin{array}{ccccccc}
1 & 2 & \dots & k & k+1 & \dots & n \\
a_{1} & a_{2} & \dots & a_{k} & b_{1} & \dots & b_{n-k} \end{array} \right)$$
with $b_{1}, \dots, b_{n-k}$ being numbers from $I_{n} \setminus S$ in
increasing order.
\end{lem}

\dowod (Theorem \ref{maintred}).
Denote by $M$ the matrix corresponding to the problem $(\grF, B)$,
by $M'$ the enlarged matrix corresponding to the problem $(\grF', B')$.
Let $c = \# P$, $s = \# B$.
Observe that $M'$ is of the following form:
$$M' = \left[ \begin{array}{cc}
M & K_{1} \\
K_{2} & N \end{array} \right],$$
where $K_{1}$ is the matrix with $s$ rows corresponding to the
conditions from $C_{\grF}$ and $c$ columns assigned to new monomials from $P$,
$K_{2}$ is the matrix with $c$ rows corresponding to the conditions
from $F$ and $s$ columns corresponding to the monomials from $B$.
The matrix $N$ is just the matrix of the $(\{F\},P)$ problem.
In the last $c$ rows we have the new indeterminates (adding $F$ to $\grF$
is adding a new independent node to interpolation).
From the generalized Laplace rule the determinant of $M'$ is the sum
of all possible $c$-minors (i. e. minors of rank $c$) from the matrix $[K_{2} \ N]$ multiplicated
by a suitable minor from the matrix $[M \ K_{1}]$
(with coefficients $1$ or $-1$) (see Lemma \ref{leibnitzlem}).
It is easy to see that every $c$-minors from $[K_{2} \ N]$ is a monomial
with coefficient (possibly equal to $0$). This monomial is determined only by
choosing $c$ columns (that is, $c$ monomials from $B'$) and is equal to the 
product of chosen $c$ monomials divided
by some monomial depending only on $F$.
The determinant of $N$ is a monomial with nonzero coefficient
(this follows from the assumption). If any other minor
gives the same monomial then the product of $c$ chosen monomials is equal
to the product of $c$ last monomials, hence this minor is zero (this
follows from the
assumption that $P$ is exceptional). Consequently considering
$\det M$ as a polynomial of new indeterminates with coefficients being
old indeterminates, the monomial
$\det N$ has a coefficient $\det M$ which is nonzero. $\hfill \square$
\remarkspace

\dowod (Lemma \ref{leibnitzlem}).
Let $S_{j}$ denote the group of permutations of $I_{j}$.
We have the correspondence
$$S_{k} \times S_{n-k} \times \S \ni
(\eta_{1}, \eta_{2}, S) \longmapsto
\overline{\eta_{1}} \circ \overline{\eta_{2}} \circ s(S) \in S_{n},$$
where
\begin{align*}
\overline{\eta_{1}}(i) = & \left\{ \begin{array}{ll}
\eta_{1}(i), & i \leq k \\
i, & i > k \end{array} \right. , \\
\overline{\eta_{2}}(i) = & \left\{ \begin{array}{ll}
\eta_{2}(i-k) + k, & i > k \\
i, & i \leq k \end{array} \right. . 
\end{align*}
Observe that the above correspondence is one-to-one.
Let $M = \big[ a_{i,j} \big]$. From the definition
$$ \det M = \sum\limits_{\sigma \in S_{n}} \sgn(\sigma) a_{1,\sigma(1)} \cdot
\cdots \cdot a_{n,\sigma(n)}.$$
We can identify $\eta_{1}$ with $\overline{\eta_{1}}$,
$\eta_{2}$ with $\overline{\eta_{2}}$. Proceeding with the deteminant
we have
$$\begin{array}{rcl} \det M & = &
\sum\limits_{S \in \S} \sum\limits_{\eta_{1} \in S_{k}, \eta_{2} \in S_{n-k}}
\sgn(s(S)) \cdot \sgn(\eta_{1}) \cdot \sgn(\eta_{2}) \cdot \\
& & \cdot a_{1,\eta_{1}(S(1))} \cdot \cdots \cdot a_{k,\eta_{1}(S(k))}
\ \cdot \ a_{k+1,\eta_{2}(S(k+1))} \cdot \cdots \cdot
a_{n,\eta_{2}(S(n))}. \end{array}$$
Now we can sum up this equation with respect to $\eta_{1}$ obtaining
the minor of $M_{1}$ given by $S$, and then we can do the same for
$\eta_{2}$ finishing the proof. $\hfill \square$
\remarkspace

Define $F_{d}^{n} = \{ \alpha \in \ndon : |\alpha| < d \}$.
The cardinality of $F_{d}^{n}$ is equal to ${n+d-1}\choose{n}$.
We need the following lemma:

\begin{lem}
\label{onepointlemma}
Let $n, d \in \N$, $n \geq 1$, $d \geq 1$. Consider the set $B$
of ${n+d-1}\choose{n}$ monomials. The exponents of monomials
from $B$ form a set $B'$ in $\ndon$. Then
the interpolation problem $(\{F_{d}^{n}\}, B)$ is generically
correct if and only if the set $B'$ does not lie on a hypersurface
of degree $d-1$. In particular, for $d=1$ $(\{F_{1}^{n}\},\{x^\alpha\})$
is generically correct.
\end{lem}

\dowod
The determinant of the matrix $M$ assigned to the problem $(\{F_{d}^{n}\},B)$
is a monomial with coefficient. It is enough to calculate this coefficient.
Let $B' = \{ a_{i} = (a_{1}^{i}, \dots, a_{n}^{i}), i=1, \dots, r \}$.
For every condition $\alpha \in F_{d}^{n}$ and point $a_{i}$ the assigned
entry in $M$ is equal to
$$\underbrace{a_{1}^{i}(a_{1}^{i}-1)(a_{1}^{i}-2) \cdots}_{\alpha_{1}} \cdot
  \underbrace{a_{2}^{i}(a_{2}^{i}-1)(a_{2}^{i}-2) \cdots}_{\alpha_{2}} \cdot \cdots \cdot
  \underbrace{a_{n}^{i}(a_{n}^{i}-1)(a_{n}^{i}-2) \cdots}_{\alpha_{n}},$$
where $i$ is the index of the column. Observe that
by adding a suitable linear combination of rows assigned to all $\beta \leq \alpha$
we can obtain each entry equal to
$$(a_{1}^{i})^{\alpha_{1}} (a_{2}^{i})^{\alpha_{2}} \dots (a_{n}^{i})^{\alpha_{n}}.$$
Since we take all $\alpha \in F_{d}^{n}$, we will find in $M$ all possible products
of $a_{1}^{i}, \dots, a_{n}^{i}$
up to degree $d-1$. Now 
\begin{align*}
\det M \neq 0 \iff & \textrm{ the rows of } M \textrm{ do not satisfy
linear equation with} \\
& \textrm{ nonzero coefficient } \\
\iff & \textrm{ the points from } B' \textrm{ do not satisfy an equation of} \\
& \textrm{ degree } d-1.
\end{align*} $\hfill \square$

\section{Interpolation on the plane.}
Assume now that $n=2$, $F_{d} := F_{d}^{2}$. For any finite sequence
$a_{1}, \dots, a_{k}$, $a_{i} \in I_{i}$ we define
{\it the diagram of type $(a_{1}, \dots, a_{k})$} by
$$F = \{ \alpha \in \N^{2} : \exists_{i \in I_{k}} \ |\alpha| = i - 1, \
\alpha_{2} < a_{i} \}.$$
Observe that the diagram of type $(0)$ is the empty diagram.
For example, $F$ of type $(1,2,2)$ is equal to
$\{(0,0), (1,0), (0,1), (2,0), (1,1)\}$.
Define the type $(\overline{a}, a_{1}, \dots, a_{k}) := (1,2, \dots, a, a_{1}, \dots, a_{k})$.
For example $(\overline{3}, 4, 3, 2) = (\overline{4}, 3, 2) = (1, 2, 3, 4, 3, 2)$.
$F_{d}$ is of type $(\overline{d}) = (1,2,3, \dots, d)$.
We say that the diagram {\it has at most $k$ steps} if it is of type
$(\overline{a}, a_{1}, \dots, a_{k})$, $a \geq a_{1} \geq a_{2} \geq \dots \geq a_{k} \geq 0$.
Any diagram with at most $k$ steps is a Ferrers diagram.
We say that the diagram $F$ is $k$-diagram if
$k(k+1)/2 \, | \, \#F$ and $F$ has at most
$k$ steps.

{\bf The lowest degree problem.}
With every Ferrers diagram $F$ we can assign a set of monomials
$S_{F} = \{ x^{\alpha} : \alpha \in F\}$. Now we restrict
our studies to the following situation:

Let $d \geq 1$. We want to solve an interpolation problem for
a sequence of diagrams $F_{d}$, that is, we want to find
a ``good'' set of monomials $B$ such that the problem
$(\{F_{d}\}_{k},B)$ is generically correct.
By ``good'' we understand the set given by a Ferrers diagram
(of cardinality $c = kd(d+1)/2$)
with at most 1 step. This restriction is natural:

For the purpose of interpolation we want to use
the set of $c$ first monomials with respect
to total degree ordering.
If the number of nodes multiplied by the cardinality of
$F_{d}$ coincides with the cardinality of some $F_{\ell}$ than
we want this $F_{\ell}$ to be a ``good'' set of monomials.
It is not always so, for example one can show that
$(\{F_{2}\}_{2}, F_{3})$ is not generically correct for
interpolation. However one can expect that for the number of
nodes large enough the problem is generically correct
for good set of monomials. We will solve this problem in the
cases $d=1,\dots,13$, i. e. we will show that interpolating
values and partial derivatives up to order $12$ can be done using
polynomials with the lowest possible degree. For $d=1,2,3$ all
initial cases will be proven here, for $4 \leq d \leq 13$ a suitable
computation can be done using a computer program.

We will say that the problem $(\{F_{d}\}_{k}, F)$ is generically correct
if the problem $(\{F_{d}\}_{k}, S_{F})$ is generically correct.
For a generically correct problem
$(\{F_{d}\}_{k}, F)$ we will say that $F$ is good for interpolating in
$k$ nodes of type $F_{d}$.

\section{Reductions.}

We say that a $d$-diagram $F$ of type $(\overline{a}, a_{1}, \dots, a_{d})$,
($a_{d} > 0$),
is {\it $d$-reducible} if the following holds:
there exist $v_{1}, \dots, v_{d} \in \N$ such that
$$v_{i} = \max \{ \ell : \ell \in \{1, \dots, d\}, \ell \leq a_{i}, \ell
\neq v_{j} \textrm{ for } i < j \leq d \}.$$
If $F$ is $d$-reducible then the diagram $r(F)$ of type
$$(\overline{a}, a_{1} - v_{1}, \dots, a_{d} - v_{d})$$
will be called a {\it $d$-reduction of $F$}.

We say that a $d$-diagram $F$ of type $(\overline{a}, a_{1}, \dots, a_{d})$
is {\it proper} if $a \geq d$ and $\forall_{i=1, \dots, k-1}$
$a_{i} = a_{i+1} \impl a_{i} \geq d$. A proper $d$-diagram is {\it safely
proper} if it is of type $(\overline{a},a_{1},\dots,a_{d})$, $a \geq 2d$, $a_{d} > 0$.

\remarkspace
\remark Observe that $(\overline{3},3,3)$ is not a proper $3$-diagram, such
as $(\overline{4},1,1)$, but $(\overline{3},3,3,3)$ is.

\remarkspace
\remark Observe how we can find
a sequence $(v_{1}, \dots, v_{d})$ for reducing proper diagram.
We start from $v_{d}$
and then define all the $v_{i}$ in decreasing order. As long as
$a_{i}$ is strictly smaller than $v_{i}$ we take $v_{i} = a_{i}$.
When $a_{i} \leq d$ for $v_{i}$ we choose the greatest number between
$1$ and $d$ that has not been used before.

If $E$ is a finite set of monomials then
by $\degred(E)$ we will denote the degree of a product of monomials from $E$,
$$\degred(E) = \deg \prod\limits_{x^{\alpha} \in E} x^\alpha.$$

Now we will show useful lemmas and a proposition:

\begin{lem}
\label{properisred}
Every proper $d$-diagram is $d$-reducible.
\end{lem}

\begin{lem}
\label{properredtoproper}
If $F$ is a safely proper $d$-diagram
then the $d$-reduction of $F$ is a proper $d$-diagram.
\end{lem}

\begin{prop}
\label{exceptionalset}
Let $F$ be a proper $d$-diagram, let $F'$ be the $d$-reduction of
$F$. Let $S$ be the set of monomials
assigned to $F$, $S'$ be the set of monomials assigned to $F'$.
Then $S \setminus S'$ is exceptional in $S$ with respect to $F_{d}$.
\end{prop}

\dowod (Lemma \ref{properisred}).
Assume that $F$ of type $(\overline{a}, a_{1}, \dots, a_{d})$ is proper
but not reducible. Then there exists $i, j \in I_{d}$, $i < j$ such that
$a_{i} = v_{j}$. But then $a_{i} \leq a_{j}$ implies $a_{i} \geq d$, so it
must exist $v_{i}$ suitable for $a_{i}$, contradiction. \hfill $\square$

\remarkspace
\dowod (Lemma \ref{properredtoproper}).
It is easy to see that the $d$-reduction of a diagram with at most
$d$ steps is a diagram with at most $d$ steps. Also $\# r(F) = \# F - d(d+1)/2$,
so $r(F)$ is again a $d$-diagram. Let $F$ be a diagram of type
$(\overline{a}, a_{1}, \dots, a_{d})$, $a_{d} > 0$, $a \geq d$.
Applying the reduction
$(v_{1}, \dots, v_{n})$ to $F$ we have the following:
$$v_{i} = a_{i} \impl v_{j} = a_{j} \textrm{ for } j > i, \qquad
v_{i} < a_{i} \impl v_{1} < v_{2} < \dots < v_{i}.$$
So some of $a_{i}$ will be cancelled to $0$, all weak inequalities
reduce to strong ones and all equalities $a_{i} + 1 = a_{i+1}$ reduce to
weak inequalities. The last may happen only for $a_{i} \geq 2d$, so
if $a_{i} - v_{i} = a_{i+1} - v_{i+1}$ then $a_{i} - v_{i} \geq d$.
\hfill $\square$

\remarkspace
\dowod (Proposition \ref{exceptionalset}).
This proposition is fundamental. Together with Theorem \ref{maintred}
it allows the ``induction step''.\\
Let $E = S \setminus S'$. We want to show that the problem
$(\{F_{d}\}, E)$ is generically correct. The exponents of monomials
(considered as points in $\N^{2}$) lie on $d$ skew lines $L_{1}, \dots, L_{d}$
(with equations $y+x+k=0$ for some $k$). We may assume that
the line $L_{i}$ contains exactly $i$ points from $E$.
Suppose that there exists
a curve $C$ of degree $d-1$ containing all these points. Then the intersection
of $L_{d}$ with $C$ has at least $d$ points, and then (by Bezout's theorem,
see \cite{F})
$L_{d}$ must be a part of $C$. Inductively the equation of $C$
must contain a product of all $L_{i}$, so it has degree at least $d$, contradiction.

Now we want to satisfy the first condition.
We will show that the method of choosing $v_{i}$ used in the reduction
gives strictly maximal possible
degree ($\degred(E)$) of product of reduced monomials not lying on a curve
of degree $d-1$. In fact we will show (by induction
with respect to $k$) that
having chosen $v_{k}, \dots, v_{d}$ the procedure described
above used for $v_{1}, \dots, v_{k-1}$ gives the maximal possible
degree. 

Assume that $v_{k}, \dots, v_{d}$ have been chosen using
the above procedure, and $v_{k-1}'\neq v_{k-1}$ is given.
If $v_{k-1}' > v_{k-1}$ then $a_{k-1} \geq d$. Consider the following cases:

$v_{k-1}' > d$ means that
more than $d$ points lie on a line, so the remaining $d(d+1)/2 - v_{k-1}'
< d(d-1)/2$ points lie on a curve of degree $d-2$ and all points
lie on a curve of degree $d-1$.

$v_{k-1}' \leq d$ means that
$v_{k-1}' = v_{i}$ for some $i \geq k$. In this case
$d(d+1)/2 - i(i+1)/2$ points lie on $d-i$ lines,
two sets (each consisting of $i$) points lie on two additional lines,
and the remaining $i(i+1)/2 - 2i = (i-1)(i-2)/2 - 1$ points
lie on a curve of degree $i-3$. In conclusion all points lie on a curve
of degree $i-3+2+(d-i) = d-1$.

The case $v_{k-1}' > v_{k-1}$ has been excluded, now assume $v_{k-1}' < v_{k-1}$.
We will apply the above method to obtain the maximal degree of product
of reduced monomials.
Consider two cases.

Case 1. $a_{k-1} \geq d$. We will choose $v_{i}$ for reduction. In the upper line
we will write the original choice, in the lower that following $v_{k-1}'$:
$$\begin{array}{cccc}
v_{d} & \dots & v_{k} & v_{k-1} \\
v_{d} & \dots & v_{k} & v_{k-1}'
\end{array}.$$
Of course the next number chosen in the upper line will be
$v_{k-2}$. In the lower line we choose the maximal possible number (by
induction we know how to choose to obtain the maximal degree), which
is now $v_{k-1}$. We can follow this until the choice of $v_{k-1}'$
in the upper line is made. If this happens, we have used the same numbers
in both lines, and from now on we are choosing the same way:
$$\begin{array}{ccccccccccc}
v_{d} & \dots & v_{k} & v_{k-1} & v_{k-2} & \dots & v_{k-\ell} & v_{k-1}' &
        v_{k-\ell-2} & \dots & v_{1} \\
v_{d} & \dots & v_{k} & v_{k-1}'& v_{k-1} & \dots & v_{k-\ell+1} & v_{k-\ell} &
        v_{k-\ell-2} & \dots & v_{1}
\end{array}$$
Denote the degree of the largest reduced monomial by $s$, let $p = s - d$.
The degree of the product while reducing as in the first line
is
$$D_{1} = \sum_{i=k}^{d} v_{i}(p+i) + \sum_{i=k-\ell}^{k-1} v_{i}(p+i)
+ v_{k-1}'(p+k-\ell-1) + \sum_{i=1}^{k-\ell-2} v_{i}(p+i).$$
The same for the second line is
$$D_{2} = \sum_{i=k}^{d} v_{i}(p+i) + v_{k-1}'(p+k-1) + \sum_{i=k-\ell}^{k-1}
v_{i}(p+i-1) + \sum_{i=1}^{k-\ell-2} v_{i}(p+i).$$
$$\begin{array}{rcccl}
D_{1} - D_{2} & = & \sum\limits_{i=k-\ell}^{k-1} v_{i}(p+i-(p+i-1)) + & & \\
 & & + v_{k-1}'(p+k-\ell-1-(p+k-1)) & = & \\
 & = & \big( \sum\limits_{i=k-\ell}^{k-1} v_{i} \big) -
v_{k-1}'\ell > \big(\sum\limits_{i=k-\ell}^{k-1} v_{k-1}' \big) - v_{k-1}'\ell & = & 0,
\end{array}$$
which proves case 1.

Case 2. $a_{k-1} < d$ means $a_{k-1} = v_{k-1}$ and for some
$k-1 > i > \ell_{1}$
we have $a_{i} = v_{i}$, so we choose this $v_{i}$ in both lines.
Now in each step we choose the largest possible number, which can be
the same in both the upper and lower line for $i = \ell_{1}, \dots, \ell_{2}$.
$$\begin{array}{ccccccccccc}
v_{d} & \dots & v_{k} & v_{k-1} & v_{k-2} & \dots & v_{\ell_{1}+1} & v_{\ell_{1}} &
v_{\ell_{1} - 1} & \dots & v_{\ell_{2}} \\
v_{d} & \dots & v_{k} & v_{k-1}'& v_{k-2} & \dots & v_{\ell_{1}+1} & v_{\ell_{1}} &
v_{\ell_{1} - 1} & \dots & v_{\ell_{2}} \end{array}$$
Originally, in the upper line, we now choose $v_{\ell_{2} - 1}$.
In the lower line we can choose $v_{k-1}$, which is now the greatest
and has not been chosen yet. Then we choose with a ``shift''
until $v_{\ell_{3}} = v_{k-1}'$ in the upper line is chosen.
When this happens the same situations occurs in both lines.
$$\begin{array}{ccccccccccccc}
 \dots & v_{k} & v_{k-1} & v_{k-2} & \dots & v_{\ell_{2}} & v_{\ell_{2}-1} & v_{\ell_{2}-2} &
\dots & v_{\ell_{3}}   & v_{\ell_{3}-1} & \dots & v_{1} \\
\dots & v_{k} & v_{k-1}'& v_{k-2} & \dots & v_{\ell_{2}} & v_{k-1}        & v_{\ell_{2}-1} &
\dots & v_{\ell_{3}+1} & v_{\ell_{3}-1} & \dots & v_{1} \end{array}$$
Let $s$ denote again the greatest degree of reduced monomial, let $p=s-d$.
Then
$$\begin{array}{ccl}
D_{1} & = & \sum\limits_{i=k}^{d} v_{i}(p+i) + v_{k-1}(p+k-1) +
\sum\limits_{i=\ell_{2}}^{k-2} v_{i}(p+i) + \\
& & + \sum\limits_{i=\ell_{3}+1}^{\ell_{2}-1} v_{i}(p+i) +
v_{\ell_{3}}(p+\ell_{3}) + \sum\limits_{i=1}^{\ell_{3}-1} v_{i}(p+i),\end{array}$$
$$\begin{array}{ccl}
D_{2} & = & \sum\limits_{i=k}^{d} v_{i}(p+i) + v_{k-1}'(p+k-1)
+ \sum\limits_{i=\ell_{2}}^{k-2} v_{i}(p+i) + \\ && +
v_{k-1}(p+\ell_{2}-1) + \sum\limits_{i=\ell_{3}+1}^{\ell_{2}-1} v_{i}(p+i-1)
+ \sum\limits_{i=1}^{\ell_{3}-1} v_{i}(p+i).\end{array}$$
$$\begin{array}{rcccl} D_{1} - D_{2} & = &
v_{k-1}(p+k-1) - v_{k-1}'(p+k-1) + & & \\
 & & + v_{\ell_{3}}(p+\ell_{3})- v_{k-1}(p+\ell_{2}-1)
+ \sum\limits_{i=\ell_{3}+1}^{\ell_{2}-1} v_{i} & > & \\ & > &
 v_{k-1}(k-\ell_{2}) + v_{\ell_{3}}(\ell_{3}-k+1) +
\sum\limits_{i=\ell_{3}+1}^{\ell_{2}-1} v_{\ell_{3}} & = & \\ & = &
v_{k-1}(k-\ell_{2}) + v_{\ell_{3}}(\ell_{3}-k+1+\ell_{2}-\ell_{3}-1) & = & \\
 & = & (k-\ell_{2})(v_{k-1} - v_{\ell_{3}}) & \geq & 0. \end{array}$$
We have shown the following: Let $F$ be a reducible $d$-diagram of
type $(\overline{a}, a_{1}, \dots, a_{d})$ with
reduction $v_{1}, \dots, v_{d}$.
Choose a set $E'$ of monomials from $F$ having the three following
properties:
\begin{enumerate}
\item $\# E' = d(d+1)/2$,
\item exponents of monomials from $E'$ do not lie on a curve of degree $d-1$,
\item $\degred(E') = \degred(E)$.
\end{enumerate}
Let $\delta = \max \{ \deg x^{\alpha} : x^{\alpha} \in F \}$.
Let
$w_{i} = \# \{ x^{\alpha} : x^{\alpha} \in E, \ \deg x^{\alpha} = i \}$.
Then $w_{\delta - j} = v_{d - j}$ for $j=0, \dots, d-1$.
In conclusion, if we want to choose the set of monomials $E'$ with above
properties, we must choose $v_{d}$ monomials with maximal degree,
$v_{d-1}$ monomials with degree $\delta - 1$ and so on. Among
chosen monomials with prescribed degrees our reduction choses the monomials
with the greatest possible product of the second coordinate. This proves
that $E$ is exceptional. \hfill $\square$
\remarkspace

Now we can formulate and prove the main technical theorem.

\begin{twr}
\label{maintheorem}
Let $p, d$ be positive integers. Assume that for every proper
$d$-diagram $F$ of cardinality $pd(d+1)/2$
the following conditions are satisfied:
\begin{enumerate}
\item $F$ is safely proper,
\item the problem $(\{F_{d}\}_{p}, F)$ is generically correct.
\end{enumerate}
Then for any $k \geq p$ the diagram $F$ of cardinality
$kd(d+1)/2$ with at most 1 step is good for interpolation,
that is the problem $(\{F_{d}\}_{k}, F)$ is generically correct.
\end{twr}

\dowod
Let $F$ be a diagram of cardinality $kd(d+1)/2$ with at most 1 step.
Assume that the problem $(\{F_{d}\}_{k}, F)$ is not generically correct.
Naturally $F$ is a proper diagram and can be reduced $k-p$ times to
a proper $d$-diagram $G$ of cardinality $pd(d+1)/2$ (Lemmas
\ref{properisred} and \ref{properredtoproper}). Each reduction
produces a diagram which is not good for interpolation (Proposition
\ref{exceptionalset} and Theorem \ref{maintred}). Hence the
problem $(\{F_{d}\}_{p}, G)$ is not generically correct,
which contradicts the assumption. \hfill $\square$

\section{Main results.}

Now we solve the problem for $d=1,2,3$ by showing all initial
cases. For \mbox{$d=1$} (Lagrange interpolation)
it is enough to observe that every $1$-diagram
is $1$-reducible and eventually reduces to the diagram
of type $(1)$.

\begin{twr}
\label{cased=2}
The problem $(\{F_{2}\}_{k}, F)$ is generically
correct for $F$ with at most 1 step if and only if $k \notin \{2,5\}$.
In other words, we can interpolate values and first order derivatives
using polynomials with the lowest possible degree in any number of points
apart from the case of $2$ or $5$ points.
\end{twr}

\dowod
We can check by direct computation that all $2$-diagrams
of cardinality $18$ are good for interpolation. However we present
here another method not requiring computation of any determinant.
A $2$-diagram of cardinality 15 is one of the following:
$F_{1}= (\overline{5})$, $F_{2} = (\overline{4}, 4, 1)$,
$F_{3} = (\overline{4}, 3, 2)$.
It is easy to see that the first diagram can be obtained as
a reduction of a proper $2$-diagram of type $(\overline{5},2,1)$ only.
The $(\overline{5},2,1)$ diagram can be obtained only from $(\overline{5}, 3, 3)$ which
is not a reduction of another proper $2$-diagram. In conclusion
if we reduce a $2$-diagram with at most one step to the diagram
of cardinality $15$ we obtain either $F_{2}$ or $F_{3}$.
So it is enough to prove that these diagrams are good for interpolation.
Both $F_{2}$ and $F_{3}$ reduce to $(\overline{4}, 2)$. Then the reduction
goes as follows:
$$(\overline{4},2) \to (\overline{3},3) \to (\overline{2},2,1)
\to (\overline{2})$$
which is obviously good for interpolating in one point.
In view of Theorem \ref{maintheorem} we have proven our statement for
$k > 5$. For $k=1,3,4$ the corresponding diagrams are
$(\overline{2})$, $(\overline{3},3)$, $(\overline{4}, 2)$ and we have shown
they are good. For $k=2, 5$ we can calculate the determinant, but the
next remark will prove that case. \hfill $\square$

\remarkspace
\remark
Consider a plane curve (can be reducible) of degree $d$. It has $(d+1)(d+2)/2$
monomials.
If $d$ is not divisible by three then the last number is divisible by three,
let $p = (d+1)(d+2)/6$. If $p > 5$ (which gives $d>4$) then the problem
$(\{F_{2}\}_{p},F_{d})$ is generically correct.
Hence for a generic set of $p$ points the determinant is nonzero, so
the only solution for a set of values and derivatives equal to 0 is
a zero polynomial. This shows that a curve of degree greater than 4
(and not divisible by 3) cannot have $(d+1)(d+2)/6$ singularities in general position.
For the case $d=4$ we can take the double conic passing through general 5 points,
and for $d=2$ the double line passing through any 2 points, so determinant
of the matrix in these cases is equal to $0$.
For $d = 3n$ we have the following: the curve of degree $d$
having $((d+1)(d+2)-2)/6$ generic singularities cannot pass through
additional generic point.

\remarkspace
\remark
We can consider the problem of interpolating values in $p_{1}$
nodes and values with first order derivatives in $p_{2}$ nodes, that is
the problem $(\{\underbrace{F_{1}, \dots, F_{1}}_{p_{1}},
\underbrace{F_{2}, \dots, F_{2}}_{p_{2}}\}, F)$, where $F$ is a complementary
diagram with at most one step.
If $p_{2} \notin \{2,5\}$ then we can first $1$-reduce the diagram $F$
$p_{1}$ times to obtain a diagram with at most 1 step and then use
Theorem \ref{cased=2}. It is easy to see that the only not $2$-reducible
diagram of cardinality greater than 2 is of type $(\overline{a}, 1, 1)$.
If it is a reduction of another diagram then $a = 1$, and the last
diagram is good for interpolating values in 3 points. We have shown
that if $p_{1} \geq 3$ or $p_{2} \notin \{2,5\}$ then $F$ is good
for interpolation. In fact only $(p_{1}, p_{2}) \in \{ (0,2), (0,5) \}$
cannot be interpolated by a diagram with at most 1 step.

\begin{twr}
\label{cased=3}
The problem $(\{F_{3}\}_{k}, F)$ is generically
correct for $F$ with at most 1 step if and only if $k \notin \{2,5\}$.
\end{twr}

\dowod
Again we will consider all $3$-diagrams of cardinality 30 (that is
diagrams for interpolating in 5 points) that can be achieved as a
sequence of reductions of a $3$-diagram with at most 1 step.
Here are the list of them:
$$\begin{array}{c}
(\overline{5},5,5,5), \quad (\overline{6},3,3,3), \quad (\overline{6},4,3,2),
\quad (\overline{6},4,4,1), \\ (\overline{6},5,3,1), \quad (\overline{6},5,4),
\quad (\overline{6},6,2,1), \quad (\overline{6},6,3). \end{array} $$
In fact there are three another diagrams of cardinality $30$
with at most $3$ steps: $(\overline{6},5,2,2)$, $(\overline{7},1,1)$,
$(\overline{7},2)$. Two of them are not proper, the last can be obtained
from one of the following:
$$F_{1} = (\overline{7},3,3,2) \qquad F_{2} = (\overline{7},4,3,1)
\qquad F_{3} = (\overline{7},5,2,1).$$
$F_{1}$ cannot be obtained from a $3$-diagram, $F_{2}$ is a reduction
of $(\overline{7},5,5,4)$ which is not a reduction of a $3$-diagram,
$F_{3}$ can be produced from $(\overline{7},6,4,4)$ which again is not
a result of reduction.

The diagram $(\overline{5},5,5,5)$\footnote{
This diagram is not safely proper, but every diagram that reduces
to it is safely proper.} reduces to
$(\overline{5},4,3,2) \to (\overline{5},3)$ which is good for interpolating
in $3$ nodes (Lemma \ref{sp53}).
Another reductions:
$$\begin{array}{ccccccccc}
(\overline{6},3,3,3) & \to & (\overline{6},2,1) & \to & (\overline{5},3), \\
(\overline{6},4,3,2) & \to & (\overline{6},3) & \to & (\overline{4},4,4) &
\to & (\overline{3},3,2,1) & \to & (\overline{3}), \\
(\overline{6},4,4,1) & \to & (\overline{6},2,1), \\
(\overline{6},5,3,1) & \to & (\overline{6},3), \\
(\overline{6},5,4) & \to & (\overline{5},5,3,1) & \to & (\overline{5},3), \\
(\overline{6},6,2,1) & \to & (\overline{6},3). \end{array}$$
The last $(\overline{6},6,3)$ diagram will be done in Lemma \ref{sp63}.
This shows the correctness of an interpolation problem for at least $6$
nodes. For $k=1,3,4$ the corresponding diagrams are
$(\overline{3})$, $(\overline{5},3)$ and $(\overline{6},3)$, which
are also good for interpolation.
\hfill $\square$

\remarkspace
\remark
Again we can consider the problem of interpolating values in $p_{1}$
point, values and first order derivatives in $p_{2}$ points,
values and derivatives up to order two in $p_{3}$ points.
It it easy to see that for suitably large $p_{i}$ for some $i \in I_{3}$
the problem is generically correct. Using the above techniques
one can show that $p_{1} > 9$, $p_{2} > 5$
or $p_{3} > 5$ is enough. All exceptional triples also can be found:
$$(0,2,0) \quad (0,5,0) \quad (0,1,1) \quad (1,1,1) \quad (0,0,2) \quad (1,0,2)
\quad (2,0,2)$$
$$(3,0,2) \quad (0,1,2) \quad (0,3,2) \quad (0,1,4) \quad (1,1,4) \quad (0,0,5).$$

\begin{lem}
\label{sp53}
The problem $(\{F_{3}\}_{3}, (\overline{5},3))$ is generically
correct.
\end{lem}

\dowod
To $F = (\overline{5},3)$ we apply the reduction $(1,3,2)$ instead of
$(1,2,3)$. Let $D$ be the degree of a product of reduced monomials.
In this case $D = 25$. Let us choose $6$ monomials from $F$ such that
the degree of their product is greater or equal to $25$. Moreover,
assume that these monomials do not lie on a conic. The only possibility
to do that is to choose $2$ monomials of degree $5$, $3$ of degree
$4$ and $1$ of degree $3$ (like in our reduction $(1,3,2)$), the other
is to choose $3$ monomials of degree $5$. But now the degree of
the first coordinate in the product of chosen monomials
is at least $12$, while originally it is equal to $10$. This shows
that the set of chosen monomials is exceptional.
$$(\overline{5},3) \stackrel{(1,3,2)}{\longrightarrow}
 (\overline{3},3,2,1) \to (\overline{3}),$$
which completes the proof. \hfill $\square$

\begin{lem}
\label{sp63}
The problem $(\{F_{3}\}_{5}, (\overline{6},6,3))$ is generically correct.
\end{lem}

\dowod
Again the first reduction will be $(1,3,2)$ reduction. In this case the
same argument works, namely the degree of the first variable
in product of chosen monomials is equal to 17,
while the product of three monomials of
maximal degree gives 18. So
$$(\overline{6},6,3) \stackrel{(1,3,2)}{\longrightarrow}
(\overline{5},5,3,1)
\to (\overline{5},3)$$
and use Lemma \ref{sp53}. $\hfill \square$
\remarkspace

Now we can formulate the main theorem for the homogeneous conditions.

\begin{twr}
\label{maintred2}
Let $1 \leq d \leq 13$.
For any set of nodes of cardinality $k$ greater than six
the interpolation problem $(\{F_{d}\}_{k},F)$ is generically
correct for $F$ with at most 1 step. Additionally, if $d \leq 9$
then the problem $(\{F_{d}\}_{6},F)$ is generically correct
for $F$ with at most $1$ step.
\end{twr}

\dowod
For $d \leq 3$ the proofs were presented here. For greater value of
$d$ more complicated computations are needed. To deal with
all initial cases we used a suitable computer program. First, it produced
all proper $d$-diagrams for $k_{1}$ points. All these diagrams,
being safely proper, were then reduced
to $\ell(k_{1}, k_{2})$ safely proper $d$-diagrams for $k_{2} < k_{1}$ points
(this operation greatly
reduced the number of determinants to be computed). To that list
all reductions of diagrams with at most 1 step for $6, \dots, k_{1}-1$
points were added. Finally the program checked all determinants.
Here is the table which contains the number of cases ($\# k_{i}$ denotes
the number of proper $d$-diagrams for $k_{i}$ points):
$$\begin{array}{ccccccc}
d & k_{1} & \# k_{1} & k_{2} & \# k_{2} & \ell(k_{1}, k_{2}) \\ \hline
2 & 6 & 2 & -- & -- & 2 \\
3 & 6 & 4 & -- & -- & 4 \\
4 & 13 & 52 & 6 & 9 & 4 \\
5 & 20 & 899 & 6 & 21 & 5 \\
6 & 20 & 6471 & 6 & 60 & 7 \\
7 & 20 & 48274 & 6 & 175 & 12 \\
8 & 20 & 369629 & 6 & 524 & 23 \\
9 & 20 & 2887492 & 6 & 1571 & 41 \\
10 & 16 & 5755270 & 6 & 4811 & 189 \\
11 & 14 & 15296608 & 6 & 14918 & 714 \\
12 & 13 & 53382738 & 6 & 46707 & 2349 \\
13 & 11 & 71128794 & 6 & 147123 & 19787 \\
\end{array}$$

\remark
The same method can be used for larger values of $d$, but 
the time used for computation is deteriorating.
For $d=10, \dots, 13$ the problem $(\{F_{d}\}_{6}, F)$
for $F$ with at most one step is not generically correct.
We need at least $7$ nodes.
\remarkspace

If we want to interpolate with "mixed" conditions we can use the following
Theorem.
\begin{twr}
\label{mixed}
Let $S = \{F^{j}\}$ be a finite sequence of diagrams, $F^{j} = F_{s(j)}$.
Let $n(\ell) = \# \{ j :  s(j) = \ell \}$.
Assume that for some $d, p$ the problem $(\{F_{d}\}_{p}, F)$ is
generically correct for any $d$-diagram $F$. Let $D = \max \{\ell :
n(\ell) \neq 0 \}$.
Define $h$ as the least natural number such that
$$h \geq 2D, \qquad h(h+1) > (p-1)d(d+1).$$
Take $q \in \N$ such that
$$q > \frac{h(h-1) + 2D(h-1)}{d(d+1)}.$$
If $n(d) \geq q$ then the problem $(S, F)$ is generically correct
for a suitable $F$ with at most 1 step.
\end{twr}

\dowod
We want to reduce $F$ with all needed $e$-reductions for $e \neq d$ to $F'$,
then $d$-reduce $F'$ to one of the $d$-diagrams.
To do so, we must first know that every $e$-reduction is possible for $e \leq D$.
It is true as long as the diagram being reduced
is of type
$(\overline{a}, a_{1}, \dots, a_{D})$, $a \geq 2D$. If $a < 2D$ then
the \mbox{$D$-diagram} of type $(\overline{a}, a_{1}, \dots, a_{D})$
has at most $(2D-1)2D/2 + D(2D-1)$ points. On the other hand
our diagram has at least $qd(d+1)/2$ points, contradiction.
Now assume that $F'$ is of type $(\overline{a}, a_{1}, \dots, a_{D})$.
While $F'$ has more than $d$ steps the $d$-reduction does not change
$a$. It is enough to choose $a$ such that every diagram of type
$(\overline{a},a_{1}, \dots, a_{d})$ contain more than $(p-1)d(d+1)/2$
points. The cardinality of such diagram is at least $a(a+1)/2$.
We can see that $q$ was chosen to allow both to reduce $F$ to $F'$
and then safely $d$-reduce $F'$. \hfill $\square$

\remarkspace
\remark
For $d=1$ it is enough to take
$q > 2D(2D-1)$,
for $d = D$ taking $q \geq p$ is also enough, provided that
all $d$-diagrams for $p$ points are safely proper.
\remarkspace

Now we are able to proof the Theorem \ref{curvethm}. 

\remarkspace
\dowod
Observe that, following
the notations from Theorem \ref{mixed}, if $h = 2D$ then 
$$q > \frac{h(h-1) + 2D(h-1)}{d(d+1)} = \frac{4D(2D-1)}{d(d+1)}.$$
Otherwise, if $h > 2D$ and $D \geq 2$ then
$$\frac{h(h-1)+2D(h-1)}{d(d+1)} < \frac{Dh(h-1)}{d(d+1)} \leq
\frac{D(p-1)d(d+1)}{d(d+1)},$$
so taking $q > D(p-1)$ is enough. For $D = 1$ also $q > D(p-1)$ is enough.
Now taking $D = m+1$, $d = k+1$ and $p = 7$ we complete the proof.
\hfill $\square$

\remarkspace
\remark
Here are the exact values
of $r(m,k)$ for small $m, k$:

$$\begin{array}{cc|cccccccc}
   & k & 0 & 1 & 2 & 3 & 4 & 5 & 6 & 7 \\
 m &   &   &   &   &   &   &   &   &   \\ \hline
 0 &   & 0 & - & - & - & - & - & - & - \\
 1 &   & 0 & 5 & - & - & - & - & - & - \\
 2 &   & 3 & 5 & 5 & - & - & - & - & - \\
 3 &   & 8 & 6 & 5 & 6 & - & - & - & - \\
 4 &   & 15& 6 & 5 & 6 & 5 & - & - & - \\
 5 &   & 24& 8 & 6 & 7 & 5 & 7 & - & - \\
 6 &   & 35&11 & 6 & 7 & 6 & 7 & 6 & - \\
 7 &   & 48&16 & 8 & 7 & 6 & 7 & 6 & 7 \\
\end{array}$$

All exceptions were found
by a computer program using bounds from Theorem \ref{mixed} and
reduction methods. If the reduction fails the determinant was computed.
This, together with some more sophisticated methods\footnote{for example
one can create a list of ''good'' diagrams for small number of points
and then try to reduce to one of these diagrams}, allows to
reduce the time of computation considerably.

The values computed with Theorem \ref{mixed} are certainly not optimal. 
We can better
them by refining arguments used in the proof of Theorem \ref{mixed}
or by investigating "mixed initial cases":

\begin{twr}
\label{mixini}
Let $p, d, D$ be nonzero natural numbers. Assume that every diagram
$F$ of cardinality $pd(d+1)/2$ with at most $D$ steps
has the following properties:
\begin{enumerate}
\item $F$ is safely proper (with respect to $D$),
\item the problem $(\{F_{d}\}_{p}, F)$ is generically correct.
\end{enumerate}
Let $S$ be any sequence of Ferrers
diagrams containing at least $p$ diagrams of type $(\overline{d})$
and only diagrams of type $(\overline{k})$ for $k \leq D$.
Then the diagram $F$ of suitable cardinality
with at most 1 step is good for interpolation,
that is the problem $(S, F)$ is generically correct.
\end{twr}

\dowod
Use techniques similar to that used in proof of Theorem \ref{maintheorem}. \hfill $\square$

\remarkspace
\remark
It is enough to assume that $F$ is proper (not necessarily safely proper)
with respect to $D$. Here are the values of bounds for $r(m,k)$ obtained
from Theorem \ref{mixini}:

$$\begin{array}{cc|cccccccc}
   & k & 0 & 1 & 2 & 3 & 4 & 5 & 6 & 7 \\
 m &   &   &   &   &   &   &   &   &   \\ \hline
 0 &   & 0 & - & - & - & - & - & - & - \\
 1 &   & 2 & 5 & - & - & - & - & - & - \\
 2 &   & 9 & 5 & 5 & - & - & - & - & - \\
 3 &   & 18& 6 & 5 & 6 & - & - & - & - \\
 4 &   & 30&10 & 6 & 6 & 5 & - & - & - \\
 5 &   & 45&15 & 7 & 7 & 5 & 8 & - & - \\
 6 &   & 63&21 &10 & 7 & 6 & 7 & 7 & - \\
 7 &   & 84&28 &14 & 8 & 6 & 7 & 6 & 7 \\
\end{array}$$

\section{Constrained correctness.}

For a finite sequence $S$ of Ferrers diagrams, $S = \{ F^{i} \}_{i=1}^{k}$,
$F^{i} = F_{d_{i}}$,
define the diagram
$$F_{S} = \{ (\alpha_{1} + \cdots + \alpha_{k}, \beta) \in \N^{2} :
 (\alpha_{i}, \beta) \in F^{i}, i = 1, \dots, k \}.$$
Each level in $F_{S}$ is a sum of levels of diagrams from $S$.
The cardinality of $F_{S}$ is equal to $\sum_{i=1}^{k} d(i)(d(i) + 1)/2$,
so we can consider the problem $(S, F_{S})$.

\begin{twr}
The problem $(S, F_{S})$ is generically correct.
\end{twr}

\remark
The diagram $F_{S}$ is the minimal diagram for generic interpolation
for lexicographical ordering (see \cite{ASTW}).

\remarkspace
\dowod
Let $S' = \{ F^{i} \}_{i=1}^{k-1}$. Define $R = F_{S} \setminus F_{S'}$.
It is easy to see that $R$ is a set of multiindices
from $F_{S}$ with $d(k)$ points on the lowest level ($\N \times \{0\}$),
$d(k) - j$ points on the level $\N \times \{j\}$. The same method
as in the proof of Proposition \ref{exceptionalset} can be applied.
Namely, we choose $d(k)(d(k)+1)/2$ monomials with the lowest possible
second exponent, not lying on a curve of degree $d(k)-1$.
Among them we choose the monomials with the greatest first exponent.
The only possible choice to do that is to choose the set $R$.
Also monomials from $R$ do not lie on a curve of degree $d(k)-1$.
Hence, $R$ is exceptional in $F_{S}$ with respect to $F_{d(k)}$
and we use induction. \hfill $\square$

\end{document}